# A Radial Basis Function (RBF)-Finite Difference Method for the Simulation of Reaction-Diffusion Equations on Stationary Platelets within the Augmented Forcing Method

Varun Shankar [1]*, Grady B. Wright[2], Aaron L. Fogelson[3], and Robert M. Kirby[1]

[1] *School of Computing, Univ. of Utah, Salt Lake City, UT, USA*
[2] *Department of Mathematics, Boise State Univ., Boise, ID, USA*
[3] *Departments of Mathematics and Bioengineering, Univ. of Utah, Salt Lake City, UT, USA*

## SUMMARY

We present a computational method for solving the coupled problem of chemical transport in a fluid (blood) with binding/unbinding of the chemical to/from cellular (platelet) surfaces in contact with the fluid, and with transport of the chemical on the cellular surfaces. The overall framework is the Augmented Forcing Point Method (AFM) (*L. Yao and A.L. Fogelson, Simulations of chemical transport and reaction in a suspension of cells I: An augmented forcing point method for the stationary case, IJNMF (2012) 69, 1736-52.*) for solving fluid-phase transport in a region outside of a collection of cells suspended in the fluid. We introduce a novel Radial Basis Function-Finite Difference (RBF-FD) method to solve reaction-diffusion equations on the surface of each of a collection of 2D stationary platelets suspended in blood. Parametric RBFs are used to represent the geometry of the platelets and give accurate geometric information needed for the RBF-FD method. Symmetric Hermite-RBF interpolants are used for enforcing the boundary conditions on the fluid-phase chemical concentration, and their use removes a significant limitation of the original AFM. The efficacy of the new methods are shown through a series of numerical experiments; in particular, second order convergence for the coupled problem is demonstrated. Copyright © 0000 John Wiley & Sons, Ltd.

*Correspondence to: School of Computing, Univ. of Utah, Salt Lake City, UT, USA. E-mail: shankar@cs.utah.edu

*Prepared using fldauth.cls [Version: 2010/05/13 v2.00]*

## 1. INTRODUCTION

Consider a stationary fluid in which there are a number of suspended objects on whose surfaces chemical reactions may occur. Some of the chemical may unbind from the surface of a particular object, thus entering the fluid phase, and undergo diffusion in the fluid. This chemical may bind to the surface of the same or a different one of the suspended objects, and when bound may diffuse on the surface of the suspended object. At each point on the objects, the flux of chemical to (from) that surface should exactly balance the rate of consumption (production) of the chemical on that surface. That is, there should be no flux of the chemical *across* the surfaces of the moving objects. We wish to determine the surface density (amount/area) of the chemical bound at each point of the suspended objects' surfaces and determine the concentration of chemical at points of the fluid phase as well.

The specific situation we have in mind is intravascular blood clotting. The fluid is blood plasma, and the objects are small blood cells called platelets that normally circulate in the blood. During the clotting process, platelets can become activated, allowing them to stick to one another and to the vascular wall. We have modeled this process (following [1]) using the Immersed Boundary method [2] to determine the coupled motion of the fluid and platelets. In these calculations, the no-slip condition holds on the platelet surfaces, *i.e.,* the velocity at each point of the platelet matches that of the immediately adjacent fluid. The chemicals of interest in the current work are those involved in conveying 'activation' signals between platelets [1] and those involved in the coagulation enzyme network [3]. In brief, we seek to solve a diffusion-reaction equation for each chemical on a platelet surface, where the reactions are coupled to diffusion equations in the blood around the platelet. Appropriate boundary conditions are to be satisfied at all points of the surfaces of the platelets (and external boundaries).

There are a variety of Cartesian Grid methods that can be used to solve PDEs in the presence of irregularly-shaped objects within the domain. The widely-used Immersed Boundary (IB) method introduced by Peskin [4, 5] uses a discrete delta function to spread boundary forces from the IB surface to the fluid, and then the discrete fluid dynamic equations are solved using a regular discretization on a rectangular grid everywhere in the domain. The forcing methods introduced by Goldstein [6], Mohd-Yusof [7], and Kim *et al.* [8] follow the idea of the IB method in using forces to represent objects embedded in a flow, but calculate the forces using feedback terms or numerical corrections to approximately enforce boundary conditions. Fadlun [9] introduces direct forcing without modification of the stencil; but, in the end, he applies the forcing in an implicit way by modifying the stencil at grid points near the irregular boundaries. The Immersed Interface method of LeVeque and Li [10], the Embedded Boundary (EB) method of Johansen and Collela [11], the sharp interface method of Udaykumar and coworkers [12, 13], and the capacity function finite volume method of Calhoun and LeVeque [14] all modify the stencil at grid points near the irregular surfaces. Because of the explicit inclusion of the boundary conditions in the linear system, the methods with changed stencil often have better accuracy and stability than the direct forcing methods, while the simpler grid and uniform stencils of the latter make them easier to implement and allow use of fast solvers.

In the work that motivated this paper [15], a Cartesian grid method for the case in which there was no flow and the platelets are stationary was presented. In addition, that work dealt with the problem of solving pure reaction ODEs for the bound chemical densities on the platelet surfaces that also provided the boundary conditions on the diffusing chemical concentrations in the fluid around the platelets. This method was called the Augmented Direct Forcing Method; in this work, we will refer to it as the Augmented Forcing Method (AFM). This augmented forcing approach was inspired by a similar idea for fluid-structure interaction introduced by Colonius and Taira [16] (see also [17]). It was also shown that this formulation produces the same numerical results as the ghost cell method (which explicitly modifies the stencil near irregular points), and is computationally more efficient. However, that work only focused on the solution of simple ODEs on the platelet surfaces; furthermore, the method designed in that work placed constraints on how close platelets could get to one another relative to the background grid.

In this paper, we present a new numerical methodology for the simulation of reaction-diffusion equations on 2D stationary platelets that are suspended in blood and for simulating diffusion of chemical species in the blood based on boundary conditions derived from those reaction-diffusion equations. This methodology consists of two components: a new method to solve reaction-diffusion equations on curves (1D surfaces) using radial basis function generated finite differences (RBF-FD), and a version of the augmented forcing method (for the simulation of the fluid-phase diffusion equations) modified with symmetric Hermite RBF interpolation to enforce boundary conditions. We also utilize a previously developed alternate representation for the platelets based on parametric RBF interpolation [18]. In that work, it was found that normals computed using an RBF representation of the platelets (when the platelets are oddly shaped or deformed) are far more accurate than the normals obtained from the piecewise quadratic representation that is traditionally employed (e.g., [15]). Since the accuracy of normals can affect the enforcement of boundary conditions in the AFM, this new geometric model for our platelets is an important part of our methodology.

The remainder of this paper is organized as follows. We first give a precise statement of the problem we address, along with a description of two new models for chemistry on platelet surfaces involving reaction-diffusion equations. Then we briefly describe the RBF parametric representation of the platelets developed in [18] and describe in detail our approach to the numerical approximation of these new models, detailing our RBF-FD based approach to solving these equations. We go on to describe our modifications to the Augmented Forcing Method – specifically, our use of RBF-based symmetric Hermite interpolation to enforce boundary conditions within that method. In the results section, we explain how shape parameters for the different RBF interpolants in this work were selected, and we describe a series of experiments to test the behavior and convergence of our new methodology on separate and coupled problems of platelet-surface and fluid-phase chemistry.

## 2. PROBLEM STATEMENT

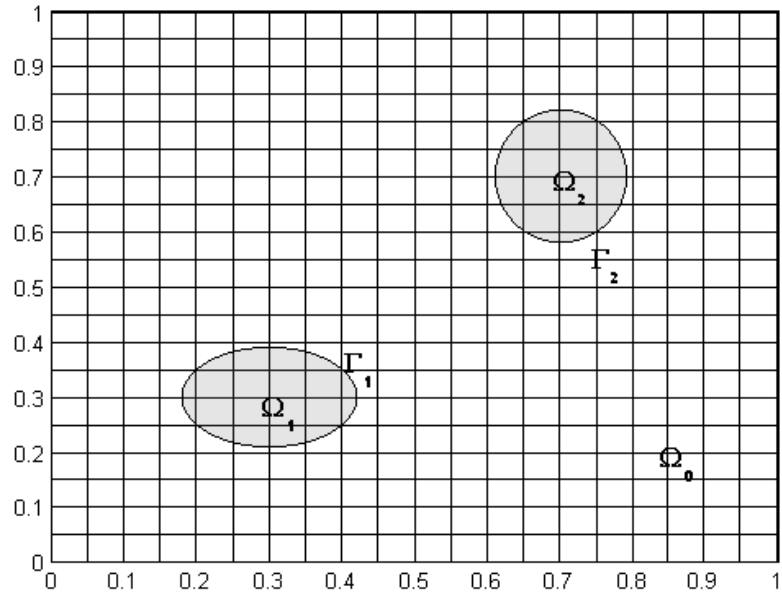

Figure 1. Illustration of a rectangular domain and grid with irregular objects embedded.

Consider a two-dimensional region $\Omega$ consisting of those points $\boldsymbol{x} = (x, y)$ within the rectangle $\Omega_0$ that are external to all of the non-overlapping subregions $\Omega_i$, $i = 1, ..., N_o$ (see Figure 1). Let $\Gamma_i$ denote the boundary of $\Omega_i$.

Let $c(\boldsymbol{x}, t) = c(x, y, t)$ be a concentration field defined for $\boldsymbol{x} \in \Omega$ and $t \geq 0$. For each $i$, let $C_i^{\mathrm{b}}(\boldsymbol{X}, t)$ and $C_i^{\mathrm{u}}(\boldsymbol{X}, t)$ be chemical density fields defined for each point $\boldsymbol{X} = (X, Y) \in \Gamma_i$ and $t \geq 0$, representing the surface densities of occupied and unoccupied binding sites respectively. We will describe these quantities in greater detail below.

We assume that $c(\boldsymbol{x}, t)$ satisfies the inhomogeneous diffusion equation

$$\frac{\partial c}{\partial t} = D\Delta c + s, \tag{1}$$

at each point $\boldsymbol{x} \in \Omega$, and that it satisfies the boundary condition

$$-D\frac{\partial c}{\partial \boldsymbol{\eta}} = k_{on}C_i^{\mathrm{u}}c - k_{off}C_i^{\mathrm{b}}. \tag{2}$$

at each point $\boldsymbol{X} \in \Gamma_i$. In these equations, $s$ is a specfied source term, $D$ is a diffusion coefficient, $k_{on}$ and $k_{off}$ are second-order and first-order rate constants, respectively, and $\boldsymbol{\eta}$ is the unit normal to $\Gamma_i$ pointing into the domain $\Omega$. Initial data for $c$ is given at all points of $\Omega$.

For the surface densities $C_i^{\mathrm{b}}$ and $C_i^{\mathrm{u}}$, we consider two variants of a reaction-diffusion model. For model 1, we imagine that the density of binding sites $C_i^{\mathrm{tot}}(\boldsymbol{X})$ is a prescribed constant at each point $\boldsymbol{X}$ on $\Gamma_i$, and we assume that the density of occupied binding sites $C_i^{\mathrm{b}}(\boldsymbol{X}, t)$, which we also refer to as the bound chemical density, satisfies

$$\frac{\partial C_i^{\mathrm{b}}}{\partial t} = k_{on}(C_i^{\mathrm{tot}} - C_i^{\mathrm{b}})c_f - k_{off}C_i^{\mathrm{b}} + D_s \Delta_{\boldsymbol{X}} C^{\mathrm{b}} b_i \tag{3}$$

at each point $\boldsymbol{X} \in \Gamma_i$. Here, $c_f$ is the value of the fluid-phase chemical concentration in the fluid adjacent to $\boldsymbol{X}$, $D_s$ is the surface diffusion coefficient, and $\Delta_{\boldsymbol{X}}$ is the Laplace-Beltrami operator on the surface.

In model 2, we instead consider a pair of coupled reaction-diffusion equations on each surface $\Gamma_i$

$$\frac{\partial C_i^{\mathrm{b}}}{\partial t} = k_{on}C_i^{\mathrm{u}}c_f - k_{off}C_i^{\mathrm{b}} + D_s^{\mathrm{b}}\Delta_{\boldsymbol{X}} C_i^{\mathrm{b}} \tag{4}$$

$$\frac{\partial C_i^{\mathrm{u}}}{\partial t} = -k_{on}C_i^{\mathrm{u}}c_f + k_{off}C_i^{\mathrm{u}} + D_s^{\mathrm{u}}\Delta_{\boldsymbol{X}} C_i^{\mathrm{u}}. \tag{5}$$

Here, the quantity $C_i^{\mathrm{u}}(\boldsymbol{X}, t)$ is the surface density of unoccupied binding sites at $\boldsymbol{X} \in \Gamma_i$ at time $t$, and $D_s^{\mathrm{u}}$ is the surface diffusion coefficient for these sites. In this variant, all binding sites diffuse on the surface, so the total density of binding sites at $\boldsymbol{X}$, $C_i^{\mathrm{tot}}(\boldsymbol{X}, t) = C_i^{\mathrm{b}}(\boldsymbol{X}, t) + C_i^{\mathrm{u}}(\boldsymbol{X}, t)$, can change in time. The setup of model 2 is intended to better represent the biological fact that the binding sites are proteins embedded in the platelet's lipid membrane and that both occupied and unoccupied proteins may diffuse. For both problems, initial data for $C_i^{\mathrm{b}}$ and $C_i^{\mathrm{u}}$ are given at all points on $\Gamma_i$.

The surface chemistry and chemical transport are coupled to that in the fluid because of the appearance of $c_f$ in the surface Equations (3) or (4) and (5), and because of the appearance of $C_i^{\mathrm{b}}$ and $C_i^{\mathrm{u}}$ in the boundary conditions (2) for the fluid-phase chemical. (For Problem 1, we set $C_i^{\mathrm{u}} = C_i^{\mathrm{tot}} - C_i^{\mathrm{b}}$ for each point $\boldsymbol{X} \in \Gamma_i$.) In our earlier work [15], only the reaction portions of these equations were considered, that is, there was no surface diffusion, and the chemical surface densities at different points were coupled only indirectly through the diffusion of fluid-phase chemicals.

## 3. GRID AND PLATELET GEOMETRY

In this section, we first introduce some terminology relevant to the rectangular Eulerian grid and to the Augmented Forcing Method. We then discuss our RBF-based parametric model for representing platelets, and the advantages it offers over the commonly-used alternatives.

### 3.1. Grid and boundary

We overlay a uniform Cartesian grid with spacing $h$ over the domain of interest, $\Omega_0$. Let $(x_i, y_j) = (ih, jh)$ denote a point of the grid. We require $c$ only within the domain $\Omega$ but nevertheless define $c_{ij}$ for all points of the mesh. Let $t_n = n\Delta t$ be the current time, where $\Delta t$ is the timestep.

To simplify the exposition, we assume that there exists a single irregular object $\Omega_1$ with boundary $\Gamma_1$. We may now classify each grid point based on its relation to the irregular object. Grid points in the domain $\Omega$ are called *fluid points*; a grid point that is covered by the object with at least one neighboring grid point not covered by the object is called a *forcing point*; finally, the grid points covered by the object that are not forcing points are called *solid points*. We also define *boundary points*, which are points on the boundary of the object whose inward normal vectors pass through forcing points; consequently, there are as many boundary points as there are forcing points. This labeling process extends to the case when multiple irregular objects exist in the domain.

### 3.2. Geometric model for platelets

We now turn our attention to the representation of platelets. In the intended application, the dynamics of the irregular objects and the fluid in which they are immersed will be described using the IB method [5, 19, 20]. In [18], we presented a parametric RBF representation of the boundaries $\Gamma_i$ of the platelets and showed that it was more accurate and less costly than the collection of techniques used in the traditional IB method (using piecewise quadratics to compute normal vectors and using finite-differences to compute forces).

For the sake of clarity and readability, we reproduce our RBF parametric model here. The model is based on explicit parametric representations of the platelets. Since our target objects are platelets, which in 2D models are nearly elliptical or circular, we use a polar parameterization. The modeling problem can be thought of as an interpolation problem. With this paradigm in mind, we first present our notation.

We represent a platelet surface at any time $t$ parametrically by

$$\boldsymbol{X}(\lambda, t) = (X(\lambda, t), Y(\lambda, t)) \tag{6}$$

where $0 \leq \lambda \leq 2\pi$ is the parametric variable and $\boldsymbol{X}(0, t) = \boldsymbol{X}(2\pi, t)$. We explicitly track a finite set of $N_d$ points $\boldsymbol{X}_1^{\mathrm{d}}(t), \ldots, \boldsymbol{X}_{N_d}^{\mathrm{d}}(t)$, which we refer to as *data sites*. Here $\boldsymbol{X}_j^{\mathrm{d}}(t) := \boldsymbol{X}(\lambda_j, t)$, $j = 1, \ldots, N_d$, and we refer to the parameter values $\lambda_1, \ldots, \lambda_{N_d}$ as the *data site nodes* (or simply *nodes*). We construct each component of $\boldsymbol{X}$ by using a smooth parametric RBF interpolant of corresponding coordinate of the data sites as discussed in detail below.

We use the interpolant and its derivatives at another set of prescribed *sample points* or *sample sites*, which correspond to $N_s$ parameter values: $\lambda_1^{\mathrm{s}}, \ldots, \lambda_{N_s}^{\mathrm{s}}$ which may be disjoint from the data site nodes. We select these $N_s$ parameter values to be equally spaced points between $0$ and $2\pi$. Below we show how to form the interpolant, and how to evaluate the interpolant at sample sites. For details on how to compute derivatives of the interpolant, see [18].

Here, we explain how to construct an RBF interpolant $X(\lambda, t)$ using the data $(\lambda_1, X_1^{\mathrm{d}}(t)), ..., (\lambda_{N_d}, X_{N_d}^{\mathrm{d}}(t))$; the construction of $Y(\lambda, t)$ component follows in a similar manner. Let $\phi(r)$ be a scalar-valued radial kernel, whose choice we discuss below. Define $X(\lambda, t)$ by

$$X(\lambda, t) = \sum_{k=1}^{N_d} \alpha_k^X \phi\left(\sqrt{2 - 2\cos(\lambda - \lambda_k)}\right).$$ (7)

Note that the square root term in Equation (7) is the Euclidean distance between the points on the unit circle whose angular coordinates are $\lambda$ and $\lambda_k$. For the geometric modeling, we use the multiquadric (MQ) radial kernel function, given explicitly by

$$\mathrm{MQ}: \quad \phi(r) = \sqrt{1 + (\varepsilon r)^2},$$ (8)

where $\varepsilon > 0$ is called the shape parameter. To have $X(\lambda, t)$ interpolate the given data, we require that the coefficients $\alpha_k^X, k = 1, ..., N_d$ satisfy the following system of equations:

$$\underbrace{\begin{bmatrix} \phi\left(r_{1,1}\right) & \cdots & \phi\left(r_{1,N_d}\right) \\ \phi\left(r_{2,1}\right) & \cdots & \phi\left(r_{2,N_d}\right) \\ \vdots & \ddots & \vdots \\ \phi\left(r_{N_d,1}\right) & \cdots & \phi\left(r_{N_d,N_d}\right) \end{bmatrix}}_{A_{rbf}} \underbrace{\begin{bmatrix} \alpha_1^X \\ \alpha_2^X \\ \vdots \\ \alpha_{N_d}^X \end{bmatrix}}_{\vec{\alpha}_{\mathrm{d}}^X} = \underbrace{\begin{bmatrix} X_1^{\mathrm{d}}(t) \\ X_2^{\mathrm{d}}(t) \\ \vdots \\ X_{N_d}^{\mathrm{d}}(t) \end{bmatrix}}_{\vec{X}_{\mathrm{d}}(t)},$$ (9)

where $r_{j,k} = \sqrt{2 - 2\cos(\lambda_j - \lambda_k)}$. Since $r_{j,k} = r_{k,j}$, the matrix $A_{rbf}$ is symmetric. More importantly, for the MQ kernels, $A_{rbf}$ is non-singular [21, 22].

We also define an *evaluation matrix*, $B$, which when multiplied by the coefficient vectors $\vec{\alpha}_{\mathrm{d}}^X$ and $\vec{\alpha}_{\mathrm{d}}^Y$, evaluates the coordinate interpolants at the set of $N_s$ evaluation nodes, and thus defines the spatial locations of the sample sites. The entries of this matrix are given by $B_{j,k} = \sqrt{2 - 2\cos(\lambda_j^e - \lambda_k)}, j = 1, 2, \ldots, N_s, k = 1, 2, \ldots, N_d$.

We conclude this section by noting that the RBF method is preferable to other popular methods such as Fourier-based methods or piecewise quadratics for modeling platelet-like shapes, as shown in [18]. The RBF method is also more flexible in terms of parameterization of objects. For example, if one were to find that a more general ellipse provided a better parameterization of the object than a circle, then the RBF method can be naturally extended to this new parameterization. The only change to Equation (7) would be to replace the distance measure in the argument of $\phi$ with the appropriate (Euclidean) distance measure on the target object for the parametrization. For an explanation of how to compute quantities like normal vectors from the RBF model, see [18]. For a detailed study on the eigenvalue stability of the RBF formulation on both periodic and non-periodic domains, see [23].

# 4. NUMERICAL SOLUTION OF REACTION-DIFFUSION MODELS FOR PLATELET CHEMISTRY

In this section, we discuss our numerical method for simulating the two models presented in Section 2, and our choices for approximating the different components of those models. Before we proceed, it is useful to discuss our time-stepping scheme for the numerical solution to the coupled problem presented in Section 2. As in [15], we use a fractional step approach in which during each timestep we first update the surface densities $C_i^{\mathrm{b}}$ and $C_i^{\mathrm{u}}$ for each platelet $i$ using known values of the fluid-phase concentrations to determine the values of $c_f$ in Equation (3) or Equations (4)-(5). Then, using the new as well as older surface densities in the boundary conditions (2), we update the fluid-phase concentration by solving Equation (1). Hence, in describing how we advance each of the surface densities or fluid-phase concentration, we regard the other as known.

In order to obtain the numerical solution of the PDEs of models 1 and 2, several components are required. An approximation of the local fluid-phase chemical concentration $c_f$ must be obtained at each sample site; then, an approximation to the surface Laplacian must be computed; finally, stable and efficient time-stepping schemes must be selected to advance the solutions in time.

## 4.1. Interpolating fluid-phase concentrations

In [15], Moving Least Squares (MLS) was used in order to construct a smooth approximation to $c_f$ using chemical concentrations from nearby patches of fluid points. This performed better than an alternate approach with bivariate quadratic interpolation, which produced undesirable spatial oscillations in $C^{\mathrm{b}}$. However, there are two potential issues with the use of MLS. First, it requires the solution of several (small) linear systems [24], in this case constructed from the background Eulerian grid; on a very fine grid, this cost may not be trivial. Second, if two platelets are very close to one another, an insufficient number of fluid points may be available for the construction of the MLS approximation to $c_f$ for each of those platelets.

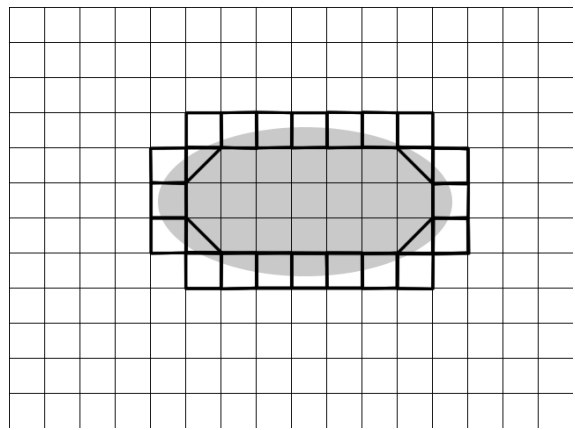

Figure 2. Illustration of the bilinear interpolation stencil for a platelet.

In order to avoid these potential issues, we use bilinear interpolation, as shown in Fig. 2; our reasoning is that bilinear interpolation has fewer degrees of freedom than bivariate quadratic interpolation while not requiring the solution of several linear systems, like MLS does. Our original approach, which we modify somewhat below, is as follows: For each sample site, we first find the Eulerian grid cell in the sample site lies. If all four corners of the cell are fluid points or forcing points, we use bilinear interpolation of the physically meaningful concentrations at the corners to the sample site. It is possible that one of the corner points is a solid point for which there is no physically meaningful concentration. In that case, we linearly interpolate concentrations from the other three corner points to this fourth corner point.

When we used this approach, our experiments showed that the resulting interpolated chemical concentration field $c_f$ was insufficiently smooth and that the overall accuracy was lower than we expected. Therefore, we instead use this procedure to first interpolate grid concentrations to *data sites* rather than sample sites. We then construct a parametric RBF interpolant of this data in the manner described in Section 3. Lastly, we evaluate the RBF interpolant at the sample sites, but with a shape parameter value slightly smaller (by a factor of 0.99), than the value used to construct the interpolant. The theory (and rationale) behind the procedure is described in detail in [25]; essentially, this procedure is a smoothing operation where we interpolate the data with a basis function, and then replace the basis function with a smoother basis function during evaluation. Figure 1 in [25] shows this process with the Multiquadric RBF, where increasing the parameter 'c' (equivalent to reducing the shape parameter $\epsilon_{geom}$ in our work to some $\epsilon_{eval}$) smoothes a noisy Lidar scan. Section 4 in [25] explains why this is equivalent to using a low-pass filter on the interpolated data by writing out the procedure in terms of convolutions against a smoother basis function. This procedure therefore gives us a smoother concentration field $c_f$ at the sample sites. In our platelet applications, we use the same parametrization for the platelets even when they are moving. This means that we can precompute the RBF interpolation and evaluation operators and use them with a single matrix-vector multiplication per platelet for each fluid-phase chemical species.

This interpolation method can be used even when two platelets are a single grid cell apart, which is an improvement over the method presented in [15] that requires that platelets be no closer than two grid widths. This feature is required for physically-relevant modeling of aggregation; in platelet aggregation simulations, platelets may indeed be very close to one another when in an aggregate.

### 4.2. Approximating the surface Laplacian

Recently RBFs have been used to compute an approximation to the surface Laplacian in the context of a pseudospectral method for reaction-diffusion equations on manifolds [26]. In that study, *global* RBF interpolants were used to approximate the surface Laplacian at a set of "scattered" nodes on a given surface. To develop a less costly method that is still sufficiently accurate for our purposes, we here use finite difference (FD)-style approximations based on RBFs for the surface Laplacian. These FD formulas are generated from RBF interpolation over *local* sets of nodes on the surface. This type of method is generally referred to as the RBF-FD method and is conceptually similar to the standard FD method with the exception that the differentiation weights enforce the exact reproduction of derivatives of shifts of RBFs (rather than derivatives of polynomials as is the case with the standard

FD method) on each local set of nodes being considered. This results in sparse matrices like in the standard FD method, but with the added advantage that the RBF-FD method can naturally handle irregular geometries. Thus, this is a technique well suited for simulating reaction-diffusion equations on platelet surfaces within a platelet aggregation simulation, where they often deform significantly from their initial discoid shapes. We note that the RBF-FD method has proven successful for a number of other applications in planar domains in two and higher dimensions (*e.g.* [27–32]), and on the surface of a sphere [33, 34], but that this is the first application of the method to more general 1D surfaces (curves).

We elect to use Cartesian coordinates rather than surface-based coordinates to formulate the surface Laplacian. This is not terribly important for 1D surfaces, but very important for generalizing our method to 2D surfaces in the future since a Cartesian-based formulation completely avoids singularities that are associated with surface-based coordinates (*e.g.* the pole singularity in spherical coordinates). Additionally, the Cartesian-based formulation is quite suitable in the context of approximation with RBFs [26]. Central to this formulation is the projection operator that takes an arbitrary 2D vector field at a point $\boldsymbol{X}$ on the surface and projects it onto the tangent line to the surface at $\boldsymbol{X}$. Letting $\boldsymbol{\eta} = (\eta_x, \eta_y)$ denote the *unit* normal vector to the surface at $\boldsymbol{X}$ (which in our applications is obtained from the RBF parametric model of the platelets described in the previous section), this operator is given by

$$\mathcal{P} = \mathcal{I} - \boldsymbol{\eta}\boldsymbol{\eta}^T = \begin{bmatrix} (1 - \eta_x \eta_x) & -\eta_x \eta_y \\ -\eta_x \eta_y & (1 - \eta_y \eta_y) \end{bmatrix} = \begin{bmatrix} \mathbf{p}_x^T \\ \mathbf{p}_y^T \end{bmatrix}, \tag{10}$$

where $\mathcal{I}$ is the 2-by-2 identity matrix, and $\mathbf{p}_x$ and $\mathbf{p}_y$ represent the projection operators in the $x$ and $y$ directions, respectively. We can combine $\mathcal{P}$ with the standard gradient operator in $\mathbb{R}^2$, $\nabla = \begin{bmatrix} \partial_x & \partial_y \end{bmatrix}^T$, to define the *surface gradient operator* $\nabla_{\boldsymbol{X}}$ in Cartesian coordinates as

$$\nabla_{\boldsymbol{X}} := \mathcal{P}\nabla = \begin{bmatrix} \mathbf{p}_x \cdot \nabla \\ \mathbf{p}_y \cdot \nabla \end{bmatrix} = \begin{bmatrix} \mathcal{G}_x \\ \mathcal{G}_y \end{bmatrix}. \tag{11}$$

Noting that $\Delta_{\boldsymbol{X}} = \nabla_{\boldsymbol{X}} \cdot \nabla_{\boldsymbol{X}}$, the surface Laplacian can then be written in Cartesian coordinates as

$$\Delta_{\boldsymbol{X}} := (\mathcal{P}\nabla \cdot)\mathcal{P}\nabla = \mathcal{G}_x \mathcal{G}_x + \mathcal{G}_y \mathcal{G}_y. \tag{12}$$

The approach we use to approximate the surface Laplacian mimics the formulation given in (12) and is conceptually similar to that based on global RBFs given in [26]. It is worth noting at this point that though the normal vector is obtained from the parametric representation of the platelet, one could certainly use normal vectors derived from level set representations or, more generally, signed-distance representations of the data sites. Since the literature on computing normal vectors on point clouds is fairly rich, we focus our attention on the exposition of the RBF-FD method, assuming that we are given reasonably smooth normal vectors. The RBF-FD method we now describe is therefore a Cartesian method that can easily handle non-parametrizeable geometries.

Given a set of $N$ nodes, we first construct discrete approximations to $\mathcal{G}_x$ and $\mathcal{G}_y$ using $n$-node RBF-FD formulas (as explained below). Letting $G_x$ and $G_y$ denote these respective discrete approximations (or differentiation matrices), we then obtain the discrete approximation to the surface Laplacian, $L$, using the matrix-products as follows:

$$L := G_x G_x + G_y G_y.$$

This approach avoids the need to compute derivatives of the normal vectors of the surfaces, but does have the effect of doubling the bandwidth of the $L$ compared to $G_x$ and $G_y$.

We explain the RBF-FD method for approximating the $\mathcal{G}_x$ component of the surface gradient in (11) as the procedure for $\mathcal{G}_y$ is similar. Without loss of generality, let the sample site where we wish to approximate $\mathcal{G}_x$ be $\boldsymbol{X}_1^s$, and let $\boldsymbol{X}_2^s, \ldots, \boldsymbol{X}_n^s$ be the $n-1$ nearest neighboring sample sites to $\boldsymbol{X}_1^s$. Given samples of a scalar valued function (say chemical density) $C(\boldsymbol{X})$ at these nodes, $C_1, \ldots, C_n$, the goal is to approximate $\mathcal{G}_x C(\boldsymbol{X})$ at $\boldsymbol{X} = \boldsymbol{X}_1^s$ using a linear combination of these samples:

$$\mathcal{G}_x C(\boldsymbol{X})\Big|_{\boldsymbol{X}=\boldsymbol{X}_1^s} \approx \sum_{i=1}^n \gamma_i C_i. \tag{13}$$

In the RBF-FD method, the differentiation weights, $\gamma_i$, are computed by enforcing that this linear combination be exact for each of the RBFs $\left\{ \phi(\|\boldsymbol{X} - \boldsymbol{X}_j^s\|) \right\}_{j=1}^n$, i.e.

$$\sum_{i=1}^n \gamma_i \phi(\|\boldsymbol{X}_i - \boldsymbol{X}_j^s\|) = \mathcal{G}_x \phi(\|\boldsymbol{X} - \boldsymbol{X}_j^s\|)\Big|_{\boldsymbol{X}=\boldsymbol{X}_1^s}, \tag{14}$$

for $j = 1, ..., n$. Note that $\|\cdot\|$ is the standard two-norm (Euclidean distance) between nodes on the surface and does not depend on any surface metrics (see [35] for a theoretical discussion on using these types of RBF approximations on general surfaces). It has also been shown through experience and studies [30, 33] that better accuracy is gained by additionally requiring that the linear combination (13) be exact for a constant. Hence, we also impose the following constraint on the weights $\gamma_i$:

$$\sum_{i=1}^n \gamma_i = \mathcal{G}_x 1\Big|_{\boldsymbol{X}=\boldsymbol{X}_1^s} = 0. \tag{15}$$

The conditions (14) and (15) can be combined into the following linear system for determining the RBF-FD weights $\gamma_i$:

$$\begin{bmatrix} \phi(\|\boldsymbol{X}_1^s - \boldsymbol{X}_1^s\|) & \cdots & \phi(\|\boldsymbol{X}_1^s - \boldsymbol{X}_n^s\|) & 1 \\ \vdots & \ddots & \vdots & \vdots \\ \phi(\|\boldsymbol{X}_n^s - \boldsymbol{X}_1^s\|) & \cdots & \phi(\|\boldsymbol{X}_n^s - \boldsymbol{X}_n^s\|) & 1 \\ 1 & \cdots & 1 & 0 \end{bmatrix} \begin{bmatrix} \gamma_1 \\ \vdots \\ \gamma_n \\ \gamma_{n+1} \end{bmatrix} = \begin{bmatrix} \mathcal{G}_x \phi(\|\boldsymbol{X} - \boldsymbol{X}_1^s\|)\big|_{\boldsymbol{X}=\boldsymbol{X}_1^s} \\ \vdots \\ \mathcal{G}_x \phi(\|\boldsymbol{X} - \boldsymbol{X}_n^s\|)\big|_{\boldsymbol{X}=\boldsymbol{X}_1^s} \\ 0 \end{bmatrix}, \tag{16}$$

where $\gamma_{n+1}$ is a dummy value that is not actually used in RBF-FD approximation after this system is solved.

The solution to (16) gives the weights for the first row of the RBF-FD differentiation matrix $G_x$ corresponding to $\boldsymbol{X}_1^s$. The weights for the second row corresponding to $\boldsymbol{X}_2^s$ are obtained by finding the $n-1$ nearest neighbors to $\boldsymbol{X}_2^s$ and solving an analogous system to (16). This procedure is repeated for determining the weights for the remaining rows of $G_x$ corresponding to sample sites $\boldsymbol{X}_3^s, \ldots, \boldsymbol{X}_{N_s}^s$. The differentiation matrix $G_y$ is obtained using the same procedure, but with the operator $\mathcal{G}_x$ in (16) replaced with $\mathcal{G}_y$. Note that each row of $G_x$ and $G_y$ contain only $n$ non-zero entries.

In all the numerical results presented in Section 6 we used the Gaussian RBF $\phi(r) = e^{-(\varepsilon r)^2}$ in (16) for computing the RBF-FD weights and set $n = 3$. In the definition of $\phi(r)$, $\varepsilon$ is again called the shape parameter. Provided it is postitive and the sample sites are distinct, the linear system (16) is guaranteed to be non-singular, which means the weights are unique. Although not presented here, we did test other RBFs (such as the multiquadric and inverse multiquadric), but found that the Gaussian generally gave better results for the experiments we ran.

### 4.3. Solving reaction-diffusion equations of Models 1 and 2

We use the discretization of the surface Laplacian just described with an implicit-explicit (IMEX) time-stepping scheme, specifically the second order accurate semi-implicit backward differentiation formula (SBDF2) [36]. For model 1, Equation (3), this corresponds to the following discretization:

$$
\begin{aligned}
\left( I - \frac{2}{3} \Delta t D_s L \right) C^{n+1} = & \frac{4}{3} \left( C^n + \Delta t k_{\mathrm{on}} \left( C^{\mathrm{tot}} - C^n \right) c_f^n - \Delta t k_{\mathrm{off}} C^n \right) \\
& - \frac{1}{3} \left( C^{n-1} + 2\Delta t k_{\mathrm{on}} \left( C^{\mathrm{tot}} - C^{n-1} \right) c_f^{n-1} - 2\Delta t k_{\mathrm{off}} C^{n-1} \right),
\end{aligned}
\tag{17}
$$

where $\Delta t$ is the time step, and $C^{n+1}$, $C^n$, and $C^{n-1}$ denote vectors containing values of the density of unoccupied surface binding sites at the $N_s$ sample sites and at time steps $n+1$, $n$, and $n-1$, respectively. Note that (17) results in an $N_s \times N_s$ sparse system of equations to solve for $C^{n+1}$. Since $N_s$ is small, we opt to use a direct method to solve this system of equations, although an iterative method such as BiCGSTAB could have also been used. We note that we bootstrap (17) with one step of SBDF1 in the initial timestep.

The discretization for Equations (4)-(5) in model 2 is similar, but contains a pair of coupled equations. However, the implicit systems that result in these two equations are not coupled, since the coupling is purely through the reaction terms, which are discretized explicitly in time.

## 5. RBF HERMITE INTERPOLATION FOR THE AUGMENTED FORCING METHOD

The Augmented Forcing Method (AFM) was developed to solve the problem of simulating chemical diffusion for stationary fluid and platelets, with the end goal of eventually simulating chemical transport for full platelet simulations within the IB method [15]. In this section, we present our modifications to the AFM based on RBF Hermite interpolation.

Hermite (or more generally Hermite-Birkhoff) interpolation refers to the interpolation of data *and* derivatives of the data. While there are a variety of methods for this task in 1D (*e.g.*, global polynomials or piecewise cubic polynomials), they often are difficult or impossible to generalize to higher dimensions, especially when the data locations are non-uniform. For these problems, the RBF Hermite interpolation method [21, Chapter 36] offers a powerful solution that can be generalized to scattered nodes and various differential and integral constraints (*e.g.*, see [30, 37]). In this work, we apply so called symmetric RBF Hermite interpolation to the problem of enforcing boundary conditions on the fluid-phase chemical concentrations within the AFM, exploiting the general nature of this formulation to overcome the separation constraints imposed by the AFM on the irregular boundaries (platelets) and similar issues that can arise in handling concavities in platelet shapes.

### 5.1. The Augmented Forcing method

The key idea of the AFM (as presented in [15]) is to solve a discrete PDE at all $N_T$ grid points, except that at forcing points, the discrete PDE is modified by the addition of a forcing term that enforces boundary conditions.

We discretize the PDE for fluid-phase chemical concentrations using the standard second-order five-point stencil for the Laplacian in space and the second-order Crank-Nicolson scheme in time. Let $A$ be the matrix formed from the discretization of Equation (1) and $\mathbf{r}$ be the right hand side vector from that discretization; let $P$ be an $N_T \times N_F$ matrix that maps each forcing point index to the index of the corresponding grid point in the overall ordering of the grid unknowns used in the vector of chemical concentrations $\mathbf{c}$, *i.e.*, all the entries of $P$ are zero except for those locations corresponding to forcing point locations, which are one. Let $\mathbf{F}$ be a vector whose $N_F$ entries contain the forcing values. Let $E$ be an $N_F \times N_T$ matrix that enforces boundary conditions as described below. Then, we require the solution of the following block system of equations:

$$\begin{pmatrix} A & P \\ E & 0 \end{pmatrix} \begin{pmatrix} \mathbf{c} \\ \mathbf{F} \end{pmatrix}^{n+1} = \begin{pmatrix} \mathbf{r} \\ \mathbf{r}_{\mathrm{bc}} \end{pmatrix}. \tag{18}$$

This system is solved in two stages.

- First, we find $\mathbf{F}$ by solving the Schur Complement system of the above block system using the BiCGSTAB iterative method. This system is as follows:

$$-EA^{-1}P\mathbf{F} = \mathbf{r}_{\mathrm{bc}} - EA^{-1}\mathbf{r}. \tag{19}$$

- Having solved for $\mathbf{F}$, we then solve for the chemical concentrations by solving $A\mathbf{c} = \mathbf{r} - P\mathbf{F}$. We use a conjugate gradient solver preconditioned by the modified incomplete Cholesky factorization of $A$.

In the AFM as implemented in [15], for each forcing point, the boundary condition at the corresponding boundary point (see Section 3) and the concentrations at five nearby fluid points are used to construct a bivariate quadratic interpolant tht satisfies the boundary condition at the boundary point. A value of $c$ at the forcing point is obtained by evaluating this polynomial at the forcing point. Since the fluid concentrations are still to be determined, this gives an implicit relationship between the forcing point concentration and those at the five fluid points. This relationship is used to populate one row of the matrix $E$. With this approach, if two platelets are close to one another or the shape of the platelet is concave, there may not be a sufficient number of points necessary to perform this interpolation. In such a case, grid refinement is necessary to introduce sufficient spacing between objects.

In the following sub-sections, we present our methods for computing the prescribed boundary conditions $\mathbf{r}_{\mathrm{bc}}$ and the matrix $E$ that enforces those boundary conditions.

### 5.2. Computing $\mathbf{r}_{\mathrm{bc}}$

Upon rearranging the boundary condition given in Equation (2), we obtain the equation

$$\left(-D\frac{\partial}{\partial \boldsymbol{\eta}} - k_{\mathrm{on}}C^{\mathrm{u}}\right)c = -k_{off}C^{\mathrm{b}}. \tag{20}$$

We use this condition at each boundary point when updating the fluid-phase concentration $c^n$ to $c^{n+1}$. Because of our fractional-step approach to time-stepping, values of $C^{\mathrm{u}}$ and $C^{\mathrm{b}}$ are known at the needed times at the locations of the sample sites (for Model 1, we set $C^{\mathrm{u}} = C^{\mathrm{tot}} - C^{\mathrm{b}}$.) These are used to compute values at the boundary points, as described below, and so we can think of them as known in Equation (20) and regard this equation as a Robin condition on $c^{n+1}$. For later reference, we define the Robin boundary condition operator by $\mathcal{D} = -D\frac{\partial}{\partial \boldsymbol{\eta}} - k_{\mathrm{on}}C^{\mathrm{u}}$.

The right hand side of Equation (20) describes the known boundary conditions on the fluid-phase chemical concentration and therefore values of it at the boundary points define the vector $\mathbf{r}_{\mathrm{bc}}$. It is important to note that $\mathbf{r}_{\mathrm{bc}}$ is modified by the procedure to compute the matrix $E$. It is this modified $\mathbf{r}_{\mathrm{bc}}$ that makes its way into the right hand side of Equation (18).

Since boundary conditions are enforced at *boundary points* in the AFM, we require require values of $C^{\mathrm{b}}$ and $C^{\mathrm{u}}$ at those points. There are many ways this can be done. For example, in [15], piecewise quadratic interpolants were fit to the concentrations at the IB points, then evaluated at the boundary points. However, in this work, we have two considerations when making this choice. First, we would like the resulting concentration field to be smooth enough to ensure that the overall convergence of our method is not affected. Second, we require that the interpolant (or more generally, approximant) be efficient to compute and evaluate. This rules out directly interpolating concentrations at the sample sites as the number of sample sites ($N_s$) is much greater than the number of data sites ($N_d$) in our geometric model; it would be more efficient to construct an approximant that had as many coefficients as the number of data sites. With these considerations in mind, we determined that a parametric least-squares fit using the RBF geometric model, described below, would be a good choice.

Let $\vec{C} = [C_1, C_2, \ldots, C_{N_s}]^T$ be a vector of function values at the sample sites, representing either $C^{\mathrm{u}}$ or $C^{\mathrm{b}}$ values at those sites. Recall from Section 3 that $B$ is an $N_s \times N_d$ RBF evaluation matrix.

When $B$ is applied to the known vector of coefficients of an RBF interpolant of some quantity defined at the data sites, we obtain values of that quantity at the sample sites. Here we use $B$ in a somewhat different way; as the coefficient matrix in a least-squares problem. We seek $N_d$ coefficients $\vec{g}$ that minimize the quantity $||B\vec{g} - \vec{C}||_2^2$; that is, we seek the coefficients of the $N_d$-term RBF expansion that best fits the $N_s$ sample-site function values that are contained in $\vec{C}$. Since the matrix $B$ depends only on the fixed parameter nodes of the data sites and sample sites, and not on their actual spatial locations, it does not change in time. Thus we precompute QR-decomposition of $B$ re-use in each timestep to solve these least squares problems.

Once we obtain the coefficients $\vec{g}$, we evaluate the least-squares approximant at the boundary points, giving us values of the chemical surface densities $C^{\mathrm{u}}$ and $C^{\mathrm{b}}$ at the $N_F$ boundary points. This is done by building an $N_F \times N_d$ evaluation matrix, $\hat{B}$, and applying it to the coefficient vector $\vec{g}$. In the current paper, since the platelets are stationary, $\hat{B}$ can be precomputed and reused every time-step. In the more general problem where platelets are advected and deformed by a background flow, $\hat{B}$ must be recomputed every time-step, since the parameter values corresponding to the boundary points change as the platelets move relative to the background grid.

### 5.3. Enforcing boundary conditions with matrix $E$

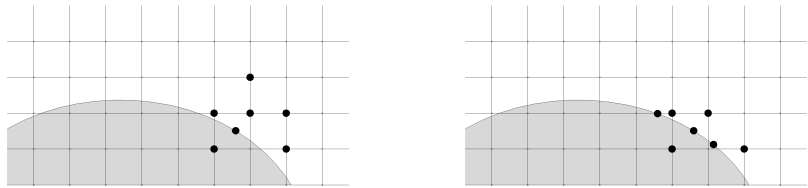

Figure 3. The figure on the left shows the number of fluid points and boundary points used in the original AFM. The figure on the right shows the number of fluid points and boundary points used within the modified AFM.

We now introduce an alternative method for computing the interpolation matrix $E$, which we refer to as $E_{\mathrm{rbf}}$. Our technique for enforcing boundary conditions is conceptually similar to the technique used in the original AFM. However, there are some significant differences, illustrated in Figure 3. For each forcing point, we now choose *three*, rather than five, nearby fluid points immediately outside the boundary to use in constructing an interpolant that satisfies the boundary conditions. Additionally, instead of using the boundary condition at a single boundary point, we now use the boundary conditions at *three* boundary points in constructing our interpolant. These changes (in conjunction with the bilinear interpolation scheme outlined in Section 3) allow platelets simulated by our modified AFM to be as close as a grid cell width apart, something which the original AFM does not allow.

Suppose we wish to impose Robin boundary conditions for $c$ along $\Gamma_1$ using the Robin boundary operator $\mathcal{D}$ from the previous sub-section. Let forcing point $\boldsymbol{B}$ have coordinates $(x_a, y_a)$ and let the corresponding boundary point $\boldsymbol{p}_a$ have coordinates $(X_a, Y_a)$. As a prototypical example, consider the layout of points in Figure 4. Here, $\boldsymbol{p}_b$ and $\boldsymbol{p}_c$ are the two boundary points closest to $\boldsymbol{p}_a$, and $\boldsymbol{p}_1$, $p_2$ and $p_3$ are the three fluid points that are neighbors to the forcing point $B$. Let $c_1$, $c_2$ and $c_3$ be the chemical concentrations at those fluid points, and recall that the boundary conditions are known at $p_a$ and $p_b$ as at $p_c$. We note that fluid points need not necessarily be chosen as shown in Figures 3 and 4; our method only requires that the selected fluid points be close to the boundary of the platelet.

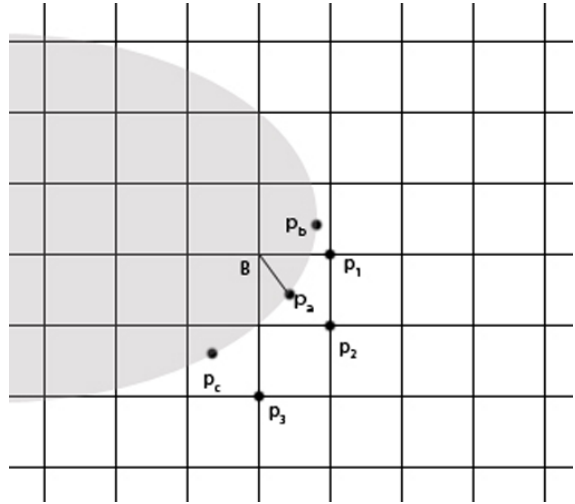

Figure 4. Illustration of the symmetric Hermite RBF interpolation stencil.

We use the symmetric RBF Hermite interpolation technique to obtain an expression for the chemical concentration at each forcing point. In this approach we construct interpolants of the form:

$$s_B(p) = \sum_{i=1}^{3} a_i \phi(\|p - p_i\|) + b_1 \mathcal{D}_{p_a} \phi(\|p - p_a\|) + b_2 \mathcal{D}_{p_b} \phi(\|p - p_b\|) + b_3 \mathcal{D}_{p_c} \phi(\|p - p_c\|), \tag{21}$$

where the boundary condition operator with subscripts is defined as

$$\mathcal{D}_{p_\iota} \phi(\|p - p_\iota\|) := \mathcal{D}\phi(\|p - x\|)\Big|_{x = p_\iota}, \quad \iota = a, b, c,$$

i.e. $\mathcal{D}$ acts on $\phi$ as a function of the subscript variable put on $\mathcal{D}$, with the other variable fixed. The interpolation conditions are given as

$$s_B(p_j) = c_j, \quad j = 1, 2, 3, \tag{22}$$

$$\mathcal{D}s_B(p)\Big|_{p = p_\iota} = \mathcal{D}c(p)\Big|_{p = p_\iota} = d_\iota, \quad \iota = a, b, c, \tag{23}$$

where the chemical concentrations $c_j$ are unknown ones from the end of the current timestep, and the boundary conditions $d_j$ are known. These interpolation conditions can be written as the following block 2-by-2 linear system of equations for determining the unknown coefficients, $a_i$ and $b_i$ in (21)

$$\underbrace{\begin{bmatrix} G & R \\ R^T & H \end{bmatrix}}_{V_B} \begin{bmatrix} \mathbf{a} \\ \mathbf{b} \end{bmatrix} = \begin{bmatrix} \mathbf{c}_B \\ \mathbf{d}_B \end{bmatrix}^{n+1}, \tag{24}$$

where $\mathbf{a}$ and $\mathbf{b}$ are vectors containing the unknown interpolation coefficients, $\mathbf{c}_B$ and $\mathbf{d}_B$ are vectors containing the respective chemical concentration (22) and boundary condition data (23) for the forcing point $B$, and the $n+1$ superscript denotes that the values are given at the next time-level. The matrix blocks in this system are defined as follows:

$$G = \begin{bmatrix} \phi(\|\boldsymbol{p}_1 - \boldsymbol{p}_1\|) & \phi(\|\boldsymbol{p}_1 - \boldsymbol{p}_2\|) & \phi(\|\boldsymbol{p}_1 - \boldsymbol{p}_3\|) \\ \phi(\|\boldsymbol{p}_2 - \boldsymbol{p}_1\|) & \phi(\|\boldsymbol{p}_2 - \boldsymbol{p}_2\|) & \phi(\|\boldsymbol{p}_2 - \boldsymbol{p}_3\|) \\ \phi(\|\boldsymbol{p}_3 - \boldsymbol{p}_1\|) & \phi(\|\boldsymbol{p}_3 - \boldsymbol{p}_2\|) & \phi(\|\boldsymbol{p}_3 - \boldsymbol{p}_3\|) \end{bmatrix}, \tag{25}$$

$$R = \begin{bmatrix} \mathcal{D}_{\boldsymbol{p}_a}\phi(\|\boldsymbol{p}_1 - \boldsymbol{p}_a\|) & \mathcal{D}_{\boldsymbol{p}_b}\phi(\|\boldsymbol{p}_1 - \boldsymbol{p}_b\|) & \mathcal{D}_{\boldsymbol{p}_c}\phi(\|\boldsymbol{p}_1 - \boldsymbol{p}_c\|) \\ \mathcal{D}_{\boldsymbol{p}_a}\phi(\|\boldsymbol{p}_2 - \boldsymbol{p}_a\|) & \mathcal{D}_{\boldsymbol{p}_b}\phi(\|\boldsymbol{p}_2 - \boldsymbol{p}_b\|) & \mathcal{D}_{\boldsymbol{p}_c}\phi(\|\boldsymbol{p}_2 - \boldsymbol{p}_c\|) \\ \mathcal{D}_{\boldsymbol{p}_a}\phi(\|\boldsymbol{p}_3 - \boldsymbol{p}_a\|) & \mathcal{D}_{\boldsymbol{p}_b}\phi(\|\boldsymbol{p}_3 - \boldsymbol{p}_b\|) & \mathcal{D}_{\boldsymbol{p}_c}\phi(\|\boldsymbol{p}_3 - \boldsymbol{p}_c\|) \end{bmatrix}, \tag{26}$$

and

$$H = \begin{bmatrix} \mathcal{D}_{\boldsymbol{p}_a}(\mathcal{D}_{\boldsymbol{p}_a}\phi(\|\boldsymbol{p}_a - \boldsymbol{p}_a\|)) & \mathcal{D}_{\boldsymbol{p}_a}(\mathcal{D}_{\boldsymbol{p}_b}\phi(\|\boldsymbol{p}_a - \boldsymbol{p}_b\|)) & \mathcal{D}_{\boldsymbol{p}_a}(\mathcal{D}_{\boldsymbol{p}_c}\phi(\|\boldsymbol{p}_a - \boldsymbol{p}_c\|)) \\ \mathcal{D}_{\boldsymbol{p}_b}(\mathcal{D}_{\boldsymbol{p}_a}\phi(\|\boldsymbol{p}_b - \boldsymbol{p}_a\|)) & \mathcal{D}_{\boldsymbol{p}_b}(\mathcal{D}_{\boldsymbol{p}_b}\phi(\|\boldsymbol{p}_b - \boldsymbol{p}_b\|)) & \mathcal{D}_{\boldsymbol{p}_b}(\mathcal{D}_{\boldsymbol{p}_c}\phi(\|\boldsymbol{p}_b - \boldsymbol{p}_c\|)) \\ \mathcal{D}_{\boldsymbol{p}_c}(\mathcal{D}_{\boldsymbol{p}_a}\phi(\|\boldsymbol{p}_c - \boldsymbol{p}_a\|)) & \mathcal{D}_{\boldsymbol{p}_c}(\mathcal{D}_{\boldsymbol{p}_b}\phi(\|\boldsymbol{p}_c - \boldsymbol{p}_b\|)) & \mathcal{D}_{\boldsymbol{p}_c}(\mathcal{D}_{\boldsymbol{p}_c}\phi(\|\boldsymbol{p}_c - \boldsymbol{p}_c\|)) \end{bmatrix}. \tag{27}$$

The matrices $G$ and $H$ are symmetric so that the composite matrix $V_B$ in (24) is symmetric. Moreover, for our choice of $\phi$ (again, the multiquadric RBF), the matrix is guaranteed to be non-singular provided the nodes $\boldsymbol{p}_1, \boldsymbol{p}_2, \boldsymbol{p}_3$, and $\boldsymbol{p}_a, \boldsymbol{p}_b, \boldsymbol{p}_c$ are distinct [21, Chapter 36]. Finally, we note that in our numerical tests, $V_B$ is also well-conditioned, thus allowing us to not only use very closely spaced fluid points from a fine grid to perform the interpolation, but also a low value for the shape parameter associated with the multiquadric RBF. The goal is to use the interpolant (21) to construct the matrix $E_{\text{rbf}}$ for enforcing boundary conditions on the interpolated chemical concentrations at the forcing points at time-level $n+1$. This matrix has dimensions $N_F \times N_T$, where $N_F$ is the number of forcing points and $N_T$ is the total number of grid points, and serves the same purpose as $E$ does in the original AFM matrix (18). The entries of $E_{\text{rbf}}$ can be obtained from (21) as follows. First, we express the interpolated chemical concentration at the forcing point $B$ as a linear combination of the chemical concentrations at the fluid grid points and the boundary conditions at the boundary points. The former are unknown as they are specified at time $n+1$, while the latter are known (see Section 5.2). The weights in this linear combination can be determined by noting that the value of the interpolant (21) at the forcing point $\boldsymbol{p} = B$ can be written as

$$s_B(B) = S_B \begin{bmatrix} \mathbf{a} \\ \mathbf{b} \end{bmatrix} = \underbrace{S_B V_B^{-1}}_{Q_B} \begin{bmatrix} \mathbf{c}_B \\ \mathbf{d}_B \end{bmatrix}^{n+1}, \tag{28}$$

where $S_{\boldsymbol{B}}$ is the row vector

$$S_{\boldsymbol{B}} = \begin{bmatrix} \phi\left(\|\boldsymbol{B} - \boldsymbol{p}_1\|\right) \\ \phi\left(\|\boldsymbol{B} - \boldsymbol{p}_2\|\right) \\ \phi\left(\|\boldsymbol{B} - \boldsymbol{p}_3\|\right) \\ \mathcal{D}_{\boldsymbol{p}_a}\phi(\|\boldsymbol{B} - \boldsymbol{p}_a\|) \\ \mathcal{D}_{\boldsymbol{p}_b}\phi(\|\boldsymbol{B} - \boldsymbol{p}_b\|) \\ \mathcal{D}_{\boldsymbol{p}_c}\phi(\|\boldsymbol{B} - \boldsymbol{p}_c\|) \end{bmatrix}^T.$$

Thus, $Q_{\boldsymbol{B}}$ contains the weights for the linear combination of chemical concentrations and boundary conditions. Letting $c_{\boldsymbol{B}} := s_{\boldsymbol{B}}(\boldsymbol{B})$ and $q_1, q_2, \ldots, q_6$ denote the entires of $Q_{\boldsymbol{B}}$, we next write this linear combination as

$$q_1 c_1 + q_2 c_2 + q_3 c_3 - c_{\boldsymbol{B}} = -q_4 d_a - q_5 d_b - q_6 d_c, \qquad (29)$$

where we have arranged the unknown values of the chemical concentration at time-level $n + 1$ on the left hand side and the known values of the boundary conditions on the right hand side.

The weights on the left hand side of (29) constitute the entries in one row of the evaluation matrix $E_{\mathrm{rbf}}$ corresponding to the forcing point $\boldsymbol{B}$. The columns for these entries correspond to the indices of the matching grid points for $\boldsymbol{B}$, $\boldsymbol{p}_1$, $\boldsymbol{p}_2$, and $\boldsymbol{p}_3$. Specifically, if $\boldsymbol{B}$ is the $k^{\mathrm{th}}$ boundary point and has lexicographic grid-index $j_1$, while $\boldsymbol{p}_1$ $\boldsymbol{p}_2$, and $\boldsymbol{p}_3$ have lexicographic indices $j_2$, $j_3$, and $j_4$, then the $k^{\mathrm{th}}$ row of $E_{\mathrm{rbf}}$ has non-zero entries

$$(E_{\mathrm{rbf}})_{k,j_1} = -1, \ (E_{\mathrm{rbf}})_{k,j_2} = -q_1, \ (E_{\mathrm{rbf}})_{k,j_3} = -q_2, \ \mathrm{and} \ (E_{\mathrm{rbf}})_{k,j_4} = -q_3.$$

Similarly, the vector of known boundary conditions $\mathbf{r}_{\mathrm{bc}}$ in (18) is populated with values from the right hand side of (29). Specifically,

$$(\mathbf{r}_{\mathrm{bc}})_k = -q_4 d_a - q_5 d_b - q_6 d_c,$$

where $d_a$, $d_b$, and $d_c$ depend on the location of the boundary point $\boldsymbol{B}$ (see Figure 4). Note that prior to modification, $(\mathbf{r}_{\mathrm{bc}})_k$ had the value $d_a$, which was in turn determined according to the method outlined in Section 5.2.

The above procedure is repeated for each of the $N_F$ forcing points $\boldsymbol{B}$. The resulting matrix $E_{\mathrm{rbf}}$ is clearly sparse, with at most four non-zero entries per row (the matrix $E$ from [15] has up to six non-zero entries per row). This procedure can be used to enforce Dirichlet and Neumann boundary conditions as well. We note that the matrix $E_{\mathrm{rbf}}$ serves the same function as the matrix $E$ mentioned in Section 5.1 and described in [15], but requires fewer fluid points in its construction. It is thus far more flexible in handling geometric features of the immersed objects.

# 6. RESULTS

In this section, we present the results of the numerical experiments performed to analyze the effects of the changes made to the AFM, as well as the results of experiments performed to analyze the new method for the coupled problems proposed in this paper. We first comment on the selection of the various shape parameters used in this work. Then, we analyze the properties of the RBF-FD discretization scheme for solving pure diffusion equations on platelet surfaces. Next, we examine the behavior of the modified AFM when using analytic boundary conditions (as opposed to deriving boundary conditions from the reaction-diffusion equations on the irregular boundaries). Having tested the convergence of the modified AFM with analytic boundary conditions, we examine the effect of varying the distance between forcing points and boundary points on the accuracy of the modified AFM. We then test the convergence of the combined method on two coupled problems, where the boundary conditions for the AFM are derived from platelet surface reaction model 1. Finally, we test the convergence of the combined method on a single coupled problem with boundary conditions for the AFM derived from platelet surface reaction model 2. Throughout this section, we compute absolute errors on the Cartesian grid and RMS (root mean squared) errors on the surfaces (as in [15]).

## 6.1. Selection of Shape Parameters

In this paper, we use RBFs in several contexts. Here, we list the shape parameters for each of those RBFs and describe the process of obtaining those shape parameters.

1. Geometry: the RBF used for the geometric modeling of the platelets is a parametric interpolant. For the selection of the shape parameter for this RBF, we follow the results obtained in [18]. For this paper, we set that shape parameter to $\varepsilon_{geom} = 0.9$.

2. Smoothing $c_f$: we use $\varepsilon_{geom}$ as the shape parameter for the parametric fit of $c_f$ at the data sites. To evaluate the RBF interpolant at sample sites and also smooth it, we use $\varepsilon_{eval} = 0.99\varepsilon_{geom}$.

3. Surface Laplace-Beltrami operator: local RBFs are used for computing the RBF-FD approximation to the surface Laplace-Beltrami operator. For these RBFs, for all tests, the shape parameter was set to $\varepsilon_{fd} = 35$. This choice was motivated, in part, by the desire to compensate for irregular point spacings on some of the perturbed objects in the tests. We note that the comparitively large value of $\varepsilon_{fd}$ is due to the partial dependence of the PDE to the Cartesian grid via the $c_f$ term, and also due to the fact that we are using the Gaussian RBF with small node spacings for the interpolation.

4. Chemical densities on platelet surfaces: we once again use the parametric model, albeit for a least-squares fit. We use the same value for the shape parameter as we do for the geometric modeling.

5. Hermite interpolant: we set the shape parameter of all the RBF Hermite interpolants (one for each forcing point) to $\varepsilon_{herm} = 5$. We found that a wide range of values could be used for $\varepsilon_{herm}$ without adversely affecting the accuracy of the AFM.

*6.2. Convergence of the RBF-FD solution of diffusion equations on a circle*

| Number of sample sites ($N_s$) | $L_2$ error | Order of convergence | $L_\infty$ error | Order of convergence |
|---|---|---|---|---|
| 50 | 2.0591e-03 | | 2.9106e-03 | |
| 100 | 5.0705e-04 | 2.02 | 7.1672e-04 | 2.02 |
| 200 | 1.2152e-04 | 2.06 | 1.7185e-04 | 2.06 |

Table I. The effect of geometric accuracy on the RBF-FD solution to the diffusion equation. The errors were measured against the exact solution at $t = 2$.

We test the RBF-FD method for solving a pure diffusion equation on the surface of a platelet. In order to test the effect of the geometric model on the solution of the diffusion equations on the irregular boundaries by the RBF-FD method, we prescribe an initial chemical density $C(\lambda, 0) = (\cos \lambda + \sin \lambda)$ on the unit circle with $0 \leq \lambda < 2\pi$. For $t \geq 0$, the function $C(\lambda, t) = e^{-t}(\cos \lambda + \sin \lambda)$ is then an exact solution to the diffusion equation on the circle when $D_s = 1$. We fix the number of data sites to $N_d = 50$ and vary the number of sample sites. To test the errors in the spatial discretization, we fix the time-step at $\Delta t = 10^{-4}$. The test was run from $t = 0$ to $t = 2$. The results of this test are shown in Table I. The results demonstrate that the RBF-FD solution to the diffusion equation on a circle exhibits second-order convergence in the sample site spacing. Similar experiments with irregularly-spaced points around the circle (results not shown) show that the convergence of the RBF-FD method gradually decreases to first order as the points become more irregularly spaced. However, the method appears tolerant to mildly uneven point spacings, both on the circle and on the test objects in Coupled Problems 2 and 3.

*6.3. Convergence of the modified AFM*

| Grid Size | $\Delta t$ | $L_2$ error | Order of convergence | $L_\infty$ error | Order of convergence |
|---|---|---|---|---|---|
| $32 \times 32$ | 0.0050 | 9.0012e-07 | | 3.3407e-06 | |
| $64 \times 64$ | 0.0025 | 2.2716e-07 | 1.99 | 8.8616e-07 | 1.92 |
| $128 \times 128$ | 0.00125 | 5.2988e-08 | 2.10 | 2.0742e-07 | 2.10 |

Table II. Results of a refinement study for the modified AFM. The errors were measured against a solution computed on a $256 \times 256$ grid as a gold standard. The errors were measured at $t = 3$.

Table II shows the results of a refinement study conducted with the modified AFM. The solution was taken to be $c(x, y, t) = \sin(\pi x) \sin(\pi y) e^{-\pi^2 t}$. The initial chemical concentration was presecribed using the values of this function at $t = 0$, and analytic boundary conditions were prescribed by applying the boundary condition operator $\mathcal{D}$ to $c$. These boundary conditions were enforced at boundary points on platelet surfaces in the modified AFM calculations, and, on the computational domain boundary, the exact solution satisfies periodic boundary conditions in the x-direction and Neumann boundary conditions in the y-direction (though our results were similar for Dirichlet and Neumann boundary conditions as well). For our refinement study, we used the Robin boundary condition operator $\mathcal{D} = -D\frac{\partial}{\partial \boldsymbol{\eta}} + 1$ with a diffusion coefficient $D = 0.1$. Two objects were embedded in the domain, a circle C1 and an ellipse E1. C1 has its center at (0.2,0.4) and a radius of 0.0995, while E1 has its center at (0.8,0.4), a semi-major axis of length 0.15 and a semi-minor axis of length 0.1. We compare the solution on grids of several sizes to a solution computed on a $256 \times 256$ grid. We also reduce the time-step by half for each progressively finer grid. The test was run from $t = 0$ to $t = 3$. The results demonstrate that the modified AFM exhibits second-order convergence in both space and time when analytic boundary conditions are prescribed.

### 6.4. Effect of location of forcing points on convergence of the modified AFM

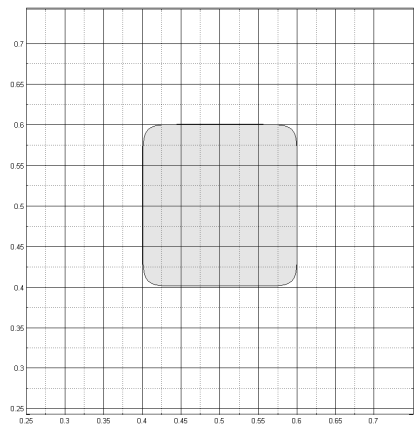

Figure 5. Illustration of the quadric object.

We wished to test whether convergence of the modified AFM is sensitive to the distance between the boundary of an irregular object and the forcing points on the grid. To accomplish this, we placed an object that looks like a square with rounded corners in the center of the domain; technically, this object is a superquadric and is shown in Fig. 5. We generated the object parametrically as follows:

$$X = x_c + r * \text{sign}(\cos \lambda)(p_x |\cos \lambda|)^m \tag{30}$$

$$Y = y_c + r * \text{sign}(\sin \lambda)(p_y |\sin \lambda|)^m \tag{31}$$

where $(x_c, y_c) = (0.5, 0.5)$, $m = 0.2$, $r = 0.0995$ and $0 \leq \lambda < 2\pi$. The test involved squeezing the sides (or top and bottom) of the object in such a way that the boundary shifts *between* grid lines without its actually crossing a grid line and thus causing generation of a new set of forcing points. We accomplished this by varying the parameters $p_x$ and $p_y$; reducing $p_x$ or $p_y$ squeezes the object either along the horizontal or the vertical, respectively, while increasing these parameters stretches the object. For this test, we successively reduced either $p_x$ or $p_y$ from 1.1 to 0.7, with $p_x = 1$ and $p_y = 1$ corresponding to the unchanged object.

We measure the error in approximating the manufactured solution to the test function $c(x, y, t) = \sin(\pi x) \sin(\pi y) e^{-\pi^2 t}$ for these different values of $p_x$ and $p_y$, and we plot the error as a function of the minimum distance between the forcing points and the boundary. We use analytic normals and sample sites locations for the rounded square so as to remove the effect of interpolation error. We set $D = 0.2$ and perform our tests on a $64 \times 64$ grid with time step $\Delta t = 0.0025$. The results are shown in Fig. 6. It is clear that the errors are unaffected by the distance between the boundary and the forcing points.

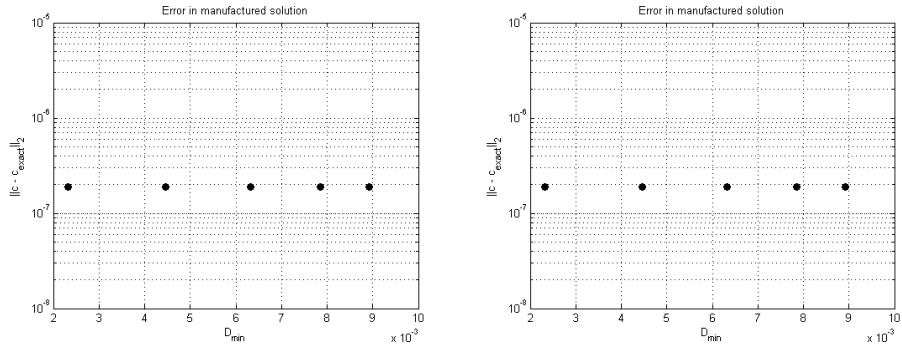

Figure 6. Error in the manufactured solution (left is for varying $p_x$ and right is for varying $p_y$) as a function of the minimum distance between a boundary point and its corresponding forcing point.

### 6.5. Convergence on coupled problems for model 1

| Grid Size | Number of sample sites ($N_s$) | $\Delta t$ | $L_2$ error | Order of convergence | $L_\infty$ error | Order of convergence |
|---|---|---|---|---|---|---|
| $32 \times 32$ | 50 | 0.0050 | 1.2841e-03 | | 1.9135e-03 | |
| $64 \times 64$ | 100 | 0.0025 | 3.2477e-04 | 1.98 | 4.9864e-04 | 1.94 |
| $128 \times 128$ | 200 | 0.00125 | 7.5756e-05 | 2.10 | 1.2041e-04 | 2.05 |

Table III. Results of a refinement study for the modified AFM on Coupled Problem 1. The errors were measured by using a solution computed on a $256 \times 256$ grid as a gold standard. The number of sample sites was also increased from $N_s = 50$ to $N_s = 200$ as the grid was refined. All errors were measured at $t = 3$.

| Number of sample sites ($N_s$) | Grid Size | $\Delta t$ | $L_2$ error | Order of convergence | $L_\infty$ error | Order of convergence |
|---|---|---|---|---|---|---|
| 50 | $32 \times 32$ | 0.0050 | 1.5567e-03 | | 2.1497e-03 | |
| 100 | $64 \times 64$ | 0.0025 | 3.6238e-04 | 2.10 | 5.0534e-04 | 2.09 |
| 200 | $128 \times 128$ | 0.00125 | 8.3943e-05 | 2.11 | 1.1706e-04 | 2.11 |

Table IV. Results of a refinement study for the RBF-FD solution to reaction-diffusion equations on the surface of platelets in Coupled Problem 1. The errors were measured using a solution computed on a $256 \times 256$ grid with analytically computed normals at $N_s = 400$ sample sites as a gold standard. The fluid grid was also refined as the number of sample sites was increased. All errors were measured at $t = 3$.

We next report on tests of the convergence of the modified AFM in conjunction with the RBF-FD method on two coupled problems. In Coupled Problem 1, the boundary conditions at boundary points for the modified AFM were obtained from the solution of reaction-diffusion equations on the surfaces of platelets C1 and E1. The diffusion coefficient for the fluid-phase chemical concentrations was set to $D = 0.1$ and that for the surface of the platelets was set to $D_s = 1$ for both platelets. The reaction rates were set to $k_{on} = 0.2$ and $k_{off} = 0.4$ for C1, and to set to $k_{on} = 0.4$ and $k_{off} = 0.2$ for E1. The fluid-phase concentrations were initialized to $c(x, y, 0) = \sin(\pi x)\sin(\pi y)$ while the platelet densities were initialized to $C(\lambda, 0) = \cos(\lambda)$, for $0 \le \lambda < 2\pi$ for both C1 and E1. The test was run from $t = 0$ to $t = 3$. Convergence was measured for both the fluid-phase concentrations and the platelet-surface concentrations. The results shown in Table III and Table IV demonstrate that the modified AFM with boundary conditions derived from the RBF-FD solution of reaction-diffusion equations on simple platelet surfaces exhibits second-order convergence in both space and time on Coupled Problem 1.

| Grid Size | Number of sample sites ($N_s$) | $\Delta t$ | $L_2$ error | Order of convergence | $L_\infty$ error | Order of convergence |
|---|---|---|---|---|---|---|
| $32 \times 32$ | 50 | 0.0050 | 9.8668e-04 | | 1.5650e-03 | |
| $64 \times 64$ | 100 | 0.0025 | 2.5374e-04 | 1.96 | 4.1418e-04 | 1.92 |
| $128 \times 128$ | 100 | 0.00125 | 5.8373e-05 | 2.12 | 9.7048e-05 | 2.09 |

Table V. Results of a refinement study for the modified AFM on Coupled Problem 2. The errors were measured by using a solution computed on a $256 \times 256$ grid as a gold standard. The number of sample sites was also increased from $N_s = 50$ to $N_s = 200$ as the grid was refined. All errors were measured at $t = 3$.

| Number of sample sites ($N_s$) | Grid Size | $\Delta t$ | $L_2$ error | Order of convergence | $L_\infty$ error | Order of convergence |
|---|---|---|---|---|---|---|
| 50 | $32 \times 32$ | 0.0050 | 1.1976e-03 | | 1.6179e-03 | |
| 100 | $64 \times 64$ | 0.0025 | 2.7351e-04 | 2.13 | 3.6505e-04 | 2.15 |
| 200 | $128 \times 128$ | 0.00125 | 6.1624e-05 | 2.15 | 8.3978e-05 | 2.13 |

Table VI. Results of a refinement study for the RBF-FD solution to reaction-diffusion equations on the surface of platelets in Coupled Problem 2. The errors were measured using a solution computed on a $256 \times 256$ grid with analytically computed normals at $N_s = 400$ sample sites as a gold standard. The grid was refined as the number of sample sites was increased. All errors were measured at $t = 3$.

Coupled Problem 2 uses the same parameters as Coupled Problem 1, but differs from that problem in solving surface reaction-diffusion equations on the ellipse E1 and on a smoothly perturbed version of ellipse E1 that we will call PE1 (Perturbed Ellipse 1). The motivation for this test was to study the behavior of the AFM on platelets which may be oddly shaped or stretched ellipses (for example, as they may be when bound to other platelets within a clot). The points on PE1 are given by

$$\boldsymbol{X}_{PE1} = \left[ 1.0 + 0.09 \exp \left( \frac{-(1 - \cos \lambda)^2}{0.1} \right) \right] \boldsymbol{X}_{E1}. \tag{32}$$

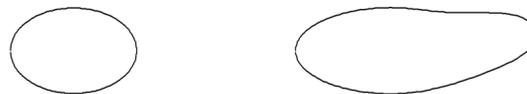

Figure 7. Illustration of the platelets in Coupled Problem 2.

The results of a convergence study of the combined method on Coupled Problem 2 are shown in Table V and Table VI. These results show that the modified AFM in conjunction with the RBF-FD method for solving reaction-diffusion equations on perturbed platelet surfaces exhibits second-order convergence in both space and time.

*6.6. Convergence on a coupled problem for model 2*

Having tested the convergence of the combined method on Coupled Problems 1 and 2 that used model 1, we now test the convergence of the combined method on a coupled problem that uses model 2. For this new coupled problem (Coupled Problem 3), we simulate the equations of model 2 on the objects E1 and PE1 using the RBF-FD method within the AFM. The reaction rates were set to the same as those in Coupled Problem 2, as were the platelet positions. The bound and unbound chemical density fields were initialized to $C^{\mathrm{b}}(\lambda) = \cos(\lambda)$ and $C^{\mathrm{u}}(\lambda) = 1 - C^{\mathrm{b}}(\lambda)$, for $0 \leq \lambda < 2\pi$, respectively. The simulation was run from $t = 0$ to $t = 3$.

| Grid Size | Number of sample sites ($N_s$) | $\Delta t$ | $L_2$ error | Order of convergence | $L_\infty$ error | Order of convergence |
|---|---|---|---|---|---|---|
| $32 \times 32$ | 50 | 0.0050 | 9.8668e-04 | | 1.5650e-03 | |
| $64 \times 64$ | 100 | 0.0025 | 2.5374e-04 | 1.96 | 4.1418e-04 | 1.92 |
| $128 \times 128$ | 200 | 0.00125 | 5.8373e-05 | 2.12 | 9.7048e-05 | 2.09 |

Table VII. Results of a refinement study for the modified AFM on Coupled Problem 3. The errors were measured by using a solution computed on a $256 \times 256$ grid as a gold standard. The number of sample sites was also increased from $N_s = 50$ to $N_s = 200$ as the grid was refined. All errors were measured at $t = 3$.

| Number of sample sites ($N_s$) | Grid Size | $\Delta t$ | $L_2$ error | Order of convergence | $L_\infty$ error | Order of convergence |
|---|---|---|---|---|---|---|
| 50 | $32 \times 32$ | 0.0050 | 1.1976e-03 | | 1.6179e-03 | |
| 100 | $64 \times 64$ | 0.0025 | 2.7351e-04 | 2.13 | 3.6505e-04 | 2.15 |
| 200 | $128 \times 128$ | 0.00125 | 6.1624e-05 | 2.15 | 8.3978e-05 | 2.13 |

Table VIII. Results of a refinement study for the RBF-FD solution to the reaction-diffusion equations for bound chemical concentrations on the surface of platelets in Coupled Problem 3. The errors were measured using a solution computed on a $256 \times 256$ grid with analytically computed normals at $N_s = 400$ sample sites as a gold standard. The grid was refined as the number of sample sites was increased. All errors were measured at $t = 3$.

| Number of sample sites ($N_s$) | Grid Size | $\Delta t$ | $L_2$ error | Order of convergence | $L_\infty$ error | Order of convergence |
|---|---|---|---|---|---|---|
| 50 | $32 \times 32$ | 0.0050 | 1.1976e-03 | | 1.6179e-03 | |
| 100 | $64 \times 64$ | 0.0025 | 2.7351e-04 | 2.13 | 3.6505e-04 | 2.15 |
| 200 | $128 \times 128$ | 0.00125 | 6.1624e-05 | 2.15 | 8.3978e-05 | 2.13 |

Table IX. Results of a refinement study for the RBF-FD solution to the reaction-diffusion equations for unbound chemical concentrations on the surface of platelets in Coupled Problem 3. The errors were measured using a solution computed on a $256 \times 256$ grid with analytically computed normals at $N_s = 400$ sample sites as a gold standard. The grid was refined as the number of sample sites was increased. All errors were measured at $t = 3$.

The results of the convergence studies are shown in Tables VII, VIII and IX. Having used the same initial conditions and platelet configurations as in Coupled Problem 2, we see identical results in terms of errors and convergence on Coupled Problem 3 for the AFM and for the PDE for $C^{\mathrm{b}}$. Furthermore, the errors and convergence for the PDE for $C^{\mathrm{u}}$ are identical to the errors and convergence for $C^{\mathrm{b}}$. We thus observe second-order convergence using our methods on Coupled Problem 3 as well. The advantage of model 2 over model 1, of course, is that one has greater flexibility in model 2, in terms of selecting initial conditions for $C^{\mathrm{u}}$ and $C^{\mathrm{b}}$, and different coefficients of diffusion as well.

## 7. DISCUSSION

The Augmented Forcing Method (AFM) was developed in [15] for the simulation of chemical transport in a stationary fluid in the presence of irregular boundaries (platelets). In that work, an ODE model for chemistry on platelet surfaces was also presented, with the ODEs contributing boundary conditions to the fluid-phase chemical diffusion equation and the fluid-phase chemical diffusion equation contributing to the ODEs. When the AFM was used in conjunction with a Crank-Nicolson timestepping method for the simulation of the combined problem, the resulting method was shown to have second-order accuracy and convergence. However, the method had the following limitations:

- the ODE model was only a simple approximation to true platelet chemistry; a reaction-diffusion PDE model would be more appropriate;
- the use of Moving Least Squares (MLS) scheme to obtain fluid-phase chemical concentrations at points on the platelet surface imposed a separation constraint on platelets – the platelets had to be at least $2h$ apart, where $h$ is the Cartesian grid spacing; and
- the AFM itself imposed another separation constraint of $2h$ on platelets because of the biquadratic interpolation stencil chosen to enforce boundary conditions on the fluid-phase chemical diffusion equation.

In this work, we introduced more complete models of platelet surface chemistry involving diffusion of chemical densities on the surface. Two models (models 1 and 2) were presented. Model 1 is a simple update to the ODE model that involved adding a surface diffusion term to the ODE (thereby giving a PDE), while model 2 aims to better match the biology of the problem by using a pair of PDEs at each point on the boundary (these PDEs are coupled to each other through their equal and opposite reaction terms).

In order to facilitate the simulation of models 1 and 2 on oddly-shaped platelets (typically seen in platelet aggregation simulations) in 2D and to remove the limitations of the AFM in its original form, we presented the following numerical methodology:

- the first application of Radial Basis Function-Finite Differences (RBF-FD) to the simulation of reaction-diffusion equations on surfaces in 2D;
- a modification to the AFM involving symmetric RBF Hermite interpolation (instead of biquadratic interpolation) to enforce boundary conditions on the fluid-phase chemical diffusion equation, thus eliminating the separation constraint on platelets simulated by the AFM; and
- a replacement for the MLS scheme used in [15] with a simple bilinear interpolation and a parametric RBF-based smoothing scheme, thereby eliminating the other separation constraint on platelets in that work.

Through numerical experimetns, we analyzed the behavior of our proposed methodology and draw the following conclusions:

- the RBF-FD approximation to the surface Laplacian, when used in conjunction with a BDF2 scheme, resulted in a method that exhibited second-order convergence in both space and time when applied to the simulation of pure diffusion equations on circles;

- the symmetric RBF Hermite interpolation scheme for enforcing boundary conditions within the AFM gave a modified AFM that also exhibited second-order convergence in both space and time for diffusion of fluid-phase concentrations; and

- the combined methodology involving RBF-FD and the AFM showed second-order convergence in both space and time on three coupled problems involving reaction-diffusion equations on platelet surfaces and a diffusion equation for the fluid-phase concentrations; Coupled Problems 1 and 2 used model 1 (simulated with SBDF2), while Coupled Problem 3 used model 2 (also simulated with SBDF2).

While we have indeed shown that the RBF-FD method can be successfully applied to the simulation of reaction-diffusion equations on platelet-like surfaces in 2D, we have yet to explore the effects on this method of using different stencil sizes, different point spacings and different RBF kernels on both 2D and 3D geometries. We plan to do this in a separate work. Also, like [15], while our results are valid for stationary platelets in stationary fluid, we have yet to explore the modified AFM and the RBF-FD method for platelets interacting with a moving fluid as simulated by the Immersed Boundary method. This, too, is the subject for future work.

**Acknowledgments:** We would like to acknowledge useful discussions concerning this work within the CLOT group at the University of Utah. The first, third and fourth authors acknowledge funding support under NIGMS grant R01-GM090203. The second author acknowledges funding support under NSF-DMS grant 1160379 and NSF-DMS grant 0934581.